\documentclass[11pt,reqno]{amsart}
\usepackage{mathrsfs, amssymb, eucal, amsmath, amsthm, amscd}
\usepackage{graphicx}

\usepackage{color}

\usepackage{bm}
\usepackage{cite}
\usepackage{hyperref}
\hypersetup{
    colorlinks=true,
    linkcolor=blue,
    filecolor=magenta,
    urlcolor=blue,
    citecolor=red,
    }
\usepackage[%
format=hang,margin=0pt,margin*=0pt,belowskip=0pt,aboveskip=0pt,indention=0pt,
parindent=0pt,hangindent=0pt,singlelinecheck=true]{caption}
\clearcaptionsetup[margin*]{singleline}
\DeclareMathOperator{\Range}{Range}

\DeclareSymbolFont{stmry}{U}{stmry}{m}{n}
\DeclareMathDelimiter\llbracket{\mathopen}{stmry}{"4A}{stmry}{"71}
\DeclareMathDelimiter\rrbracket{\mathclose}{stmry}{"4B}{stmry}{"79}

\usepackage{bbm}
\DeclareMathOperator*{\medcupd}{\mathbin{\scalebox{1.1}{\ensuremath{\bigcup}}}}%

\newcommand{\anglethick}{0.55pt}
\newcommand{\langleat}[2][1ex]{%
  \mathopen{\vcenter{\hbox{%
    \makebox[0pt][l]{\rotatebox[origin=l]{#2}{\rule{#1}{\anglethick}}}%
    \rotatebox[origin=l]{-#2}{\rule{#1}{\anglethick}}%
  }}}%
}
\newcommand{\rangleat}[2][1ex]{%
  \mathclose{\vcenter{\hbox{%
    \makebox[0pt][l]{\rotatebox[origin=l]{\numexpr180-#2\relax}{\rule{#1}{\anglethick}}}%
    \rotatebox[origin=l]{\numexpr180+#2\relax}{\rule{#1}{\anglethick}}%
  }}}%
}

\newcommand{\lmyab}{\langleat[1.2ex]{70} \mkern 0.4mu}
\newcommand{\rmyab}{\rangleat[1.2ex]{70} \mkern 1.5mu}

\allowdisplaybreaks

\numberwithin{equation}{section}
\newtheorem{theorem}{Theorem}[section]
\newtheorem{lemma}[theorem]{Lemma}
\newtheorem{proposition}[theorem]{Proposition}
\newtheorem{corollary}[theorem]{Corollary}
\newtheorem{conjecture}[theorem]{Conjecture}

\theoremstyle{remark}
\newtheorem{remark}{\bf Remark}[section]

\begin{document}
\title{Floor, ceiling and the space between}

\author[\'{A}. B\'enyi]{\'{A}rp\'{a}d B\'{e}nyi}

\address{%
Department of Mathematics, Western Washington University \\
516 High Street, Bellingham, Washington 98225, USA}

\email{Arpad.Benyi@wwu.edu}

\author[B. \'{C}urgus]{Branko \'{C}urgus}
\address{%
Department of Mathematics, Western Washington University \\
516 High Street, Bellingham, Washington 98225, USA}

\email{Branko.Curgus@wwu.edu}


\date{\today}

\begin{abstract}
Motivated by a question on the ranges of the commutators of dilated floor functions in \cite{LMR2016}, together with a related problem in \cite{XX}, we investigate the precise ranges of certain generalized polynomials dependent on a real parameter. Our analysis requires non-trivial tools, including Kronecker's approximation theorem. The results highlight sharp distinctions between irrational parameters and sub-unitary and supra-unitary rational parameters. We also propose several conjectures for the irrational and supra-unitary rational cases, supported by extensive computations in Wolfram Mathematica. 
\end{abstract}



\maketitle

\tableofcontents

\section{Introduction}

Given a real number \(x\), the \emph{floor} and \emph{ceiling} of \(x\) are defined respectively by
\[
\lfloor x\rfloor : = 
\max\bigl\{m \in \mathbb{Z} : m \leq x \bigr\}
\quad \text{and} \quad
\lceil x\rceil : = \min\bigl\{n \in \mathbb{Z} : n \geq x \bigr\}.
\]
Equivalently, \(\lfloor x\rfloor\) and \(\lceil x\rceil\) are the unique integers satisfying
\begin{equation} \label{defnflce}
\lfloor x\rfloor \leq x < \lfloor x\rfloor + 1
\quad \text{and} \quad
\lceil x\rceil - 1 < x \leq \lceil x\rceil.
\end{equation}
It is also clear from the definitions that \( \lceil x \rceil - \lfloor x \rfloor \in \{0, 1\} \), and that
\( \lceil x \rceil - \lfloor x \rfloor = 0 \) if and only if \( x \in \mathbb{Z} \).

The definitions and names of the floor and ceiling functions were introduced by Iverson in the context of computer science \cite[p. 12]{Iverson1962}; see also \cite[Chapter 3]{GKP1994} for additional history and properties. Despite their elementary definitions, these functions arise in a wide range of applications. For instance, they appear in the statement of the Quadratic Reciprocity Law, in Legendre's formula for the exponent of the highest power of a prime dividing a given integer, and in expressions involving the Riemann zeta function, among many others.

The investigation of digital straight lines, that is, drawing lines on a computer screen, naturally leads to dilated floor functions and their commutators. For $\alpha\in \mathbb R\setminus\{0\}$, define the following functions:
\begin{alignat*}{2}
\ell_{\alpha}&:\mathbb R\to\mathbb R, & \qquad \ell_\alpha (x) &= \alpha x, \\ 
g_\alpha &:\mathbb R\to\mathbb Z, & \qquad  g_\alpha (x) & =\lfloor \alpha x\rfloor. 
\end{alignat*}
From this perspective, $g_\alpha$ discretizes $\ell_\alpha$ at the length scale $1/\alpha$, in the sense that the difference $\ell_\alpha-g_\alpha$ is a bounded function. Understanding the interaction of the discretizations at two different scales $\alpha$ and $\beta$ is mathematically equivalent to the study of the commutator of the functions $g_\alpha$ and $g_\beta$ defined by
\[
[g_\alpha, g_\beta]=g_\alpha\circ g_\beta-g_\beta\circ g_\alpha.
\]
While the commutator $[\ell_\alpha, \ell_\beta]$ vanishes identically, the authors in \cite{LMR2016} highlight how  ``discretization generally destroys such commutativity,'' underscoring the new complexity introduced by the floor functions. Their main result characterizes all pairs $(\alpha, \beta)$ for which the range of the commutator $[g_\alpha, g_\beta]$ is precisely $\{0\}$; see also \cite{KMP2019}. Subsequent work in \cite{LR1,LR2} extended this analysis by asking when the range of the commutator is contained in $[0,\infty)$, drawing on tools such as Beatty sequences, Lie groups, and the Sylvester duality theorem. With the exception of \cite{LMR2016}, however, none of these studies addressed the problem of determining the \emph{precise range} of commutators.

In light of this background, it is natural to consider the simplest case where the dilation parameters coincide. Both $g_{\alpha^2}$ and $g_\alpha \circ g_\alpha$ are step-function approximations of $\ell_{\alpha^2}$. Thus, in place of the trivial commutator $[g_\alpha, g_\alpha] = 0$, one may study the difference
\[
g_{\alpha^2}-g_\alpha \circ g_\alpha,
\]
which defines a so-called generalized polynomial in the sense of \cite{BL2007}.

Specifically, for $\alpha>0$, define the function $f_\alpha:\mathbb N \to \mathbb Z$ by
\begin{equation} \label{function}
f_\alpha (n) = g_{\alpha^2}(n)-g_\alpha \bigl(g_\alpha (n) \bigr) =
\bigl\lfloor \alpha^2 \, n \bigr\rfloor - \bigl\lfloor \alpha  \lfloor \alpha n \rfloor \bigr\rfloor, \quad  \forall \, n \in \mathbb N.
\end{equation}
An instance of this function already appeared in a problem posed for high school students \cite{XX}, which asked for the precise range of $f_{\tau}$, with $\tau=\bigl(1+\sqrt{2017}\bigr)/2$.

In this note, we shall investigate a generalization of the problem in \cite{XX} and discover some nice connections with other number theoretical facts, such as the concept of modulo (or simply, mod) operator \cite{K} or Kronecker's approximation theorem \cite[Chapter III]{Cassels1957}; see also \cite{W} and \cite{Malajovich2001}. We also state several conjectures regarding the precise range of $f_\alpha$ for general supra-unitary rational parameters, and also for irrational parameters $\alpha$. Some of the questions that arise are naturally related to the size of the space between appropriate floor and ceiling functions; see, for example, Remark~\ref{rem6} and Remark ~\ref{rem23}.

\subsection{Notation} 
Above, by \(\mathbb{Z}\) we denote the set of integers and by   \(\mathbb{N}\) the set of positive integers. By \(\mathbb{R}\) we denote the set of real numbers and by \(\mathbb{Q}\) the set of rational numbers.

The \emph{fractional part} of \( x \) is defined by
\[
  \lmyab x \rmyab := x - \lfloor x \rfloor,
\]
so that \( 0 \leq \lmyab x \rmyab < 1 \) for all \( x \in \mathbb{R} \). Here we opt, as in \cite{BL2007}, for the notation \( \lmyab x \rmyab \) for the fractional part of a real number, rather than the more common \( \{x\} \), since braces are consistently used for set notation; in particular, \( \{x\} \) denotes a singleton set.

For \(a,b \in \mathbb{N}\cup \{0\}\) we set
\[
\llbracket a , b \rrbracket : = \bigl\{ x \in \mathbb{N}\cup \{0\} : a \leq x \leq b \bigr\}.
\]

We use the following operations on sets: given two sets $A, B\subseteq\mathbb R$, and given a scalar $c\in\mathbb R$, we let
\[
A+B:= \bigl\{a+b: a\in A, b\in B \bigr\},\quad \text{and}\quad cA:=\bigl\{ca: a\in A \bigr\}.
\]
We denote the set $A+(-1)B$ by $A-B$.

\section{What is the largest possible range of $f_\alpha$?}

We answer this question in the following proposition. As usual, we let
\[
\Range (f_\alpha):= \bigl\{f_\alpha (n): n\in\mathbb N \bigr\}.
\]
It is clear that for $\alpha\in\mathbb N$ we have a trivial range, $\Range (f_\alpha)=\{0\}.$ Therefore, \emph{for the remainder of our discussion we will assume that} \(\alpha > 0\) and  $\alpha\not\in\mathbb N$.

\begin{proposition} \label{largest-range}
Let $\alpha$ be a positive non-integer. Then,
\[
\Range(f_\alpha) \subseteq \bigl\llbracket 0 , \lceil\alpha\rceil \bigr\rrbracket.
\]
\end{proposition}

\begin{proof}
Write $\alpha n=\lfloor \alpha n\rfloor + \lmyab \alpha n \rmyab$. Then
\[
f_\alpha (n)
= \bigl\lfloor \alpha^2\mkern 1mu n \bigr\rfloor
- \bigl\lfloor \alpha^2 \mkern 1mu n - \alpha\mkern 1.5mu \lmyab \alpha n \rmyab \bigr\rfloor
\geq \bigl\lfloor \alpha\mkern 1.5mu \lmyab \alpha n \rmyab \bigr\rfloor
\geq 0,
\]
where we used the simple fact that, for any real numbers $x, y$, $\lfloor x\rfloor-\lfloor y\rfloor\geq\lfloor x-y\rfloor$. Furthermore, by \eqref{defnflce}, we have
\begin{align*}
f_\alpha(n)
& = \bigl\lfloor \alpha\lfloor \alpha n\rfloor + \alpha\mkern 1.5mu \lmyab \alpha n \rmyab\bigr\rfloor - \bigl\lfloor \alpha\lfloor \alpha n\rfloor \bigr\rfloor\\
& < \alpha\lfloor \alpha n\rfloor + \alpha\mkern 1.5mu \lmyab \alpha n \rmyab
-\alpha\lfloor \alpha n\rfloor + 1 \\
& = \alpha\mkern 1.5mu \lmyab \alpha n \rmyab + 1\\
& < \alpha + 1.
\end{align*}
But since $f_\alpha(n)\in\mathbb Z$, we obtain that $f_\alpha (n)\leq \lfloor \alpha +1 \rfloor=\lceil \alpha\rceil.$
\end{proof}

It is useful to note that our proof actually proves the stronger inequality
\begin{equation} \label{bounds}
\bigl\lfloor \alpha \mkern 1.5mu \lmyab \alpha n \rmyab \bigr\rfloor\leq f_\alpha (n)
\leq \bigl\lceil \alpha\mkern 1.5mu \lmyab \alpha n \rmyab\bigr\rceil, \quad \forall n\in\mathbb{N}.
\end{equation}
In particular, \eqref{bounds} implies that for all $n\in\mathbb N$ we have
\begin{equation} \label{function-revisit}
f_\alpha(n)
= \bigl\lfloor \alpha\mkern 1.5mu \lmyab \alpha n \rmyab\bigr\rfloor
\quad \text{or} \quad
f_\alpha(n) = \bigl\lceil \alpha\mkern 1.5mu \lmyab \alpha n \rmyab \bigr\rceil.
\end{equation}

\begin{corollary}
\label{01range}
If $\alpha\in (0, 1)$, then \(\{0\}\subseteq\Range(f_\alpha)\subseteq\{0, 1\}.\)
\end{corollary}

\begin{proof}
Simply note that if $\alpha\in (0, 1)$, then $f_\alpha(1)=0$ and $\lceil\alpha\rceil=1$.
\end{proof}
The previous corollary establishes that, for all sub-unitary parameters $\alpha$, we have that either $\Range(f_\alpha)=\{0\}$ or $\Range(f_\alpha)=\{0, 1\}$; see Proposition~\ref{subunitary} and Conjecture~\ref{con3}.

As we shall soon see, there is a significant difference between the mapping properties of $f_\alpha$ for sub-unitary rational parameters $\alpha$ and supra-unitary ones. Naturally, one may ask if one can reduce one case to the other, but the answer is not so simple.

The search for a better understanding of the range of $f_\alpha$ leads in a natural way to equivalent formulations of $f_\alpha$, besides the ones already expressed in \eqref{function} or \eqref{function-revisit}. When utilized carefully, these alternate formulations can provide further insight in particular situations. While we do not necessarily need these most general formulations, we will find it useful to revisit the calculations in a subsequent section dealing with rational parameters $\alpha>1$. We note also that our alternate formulations below can be obtained as a particular version of \cite[Theorem 1]{Fraenkel1994}. 

In general, writing
\[
\alpha=\lfloor\alpha\rfloor + \lmyab \alpha \rmyab,
\]
a series of routine calculations yield
\begin{equation} \label{new-f}
f_\alpha(n)
= - \lfloor \alpha\rfloor \bigl\lfloor \lmyab \alpha \rmyab\mkern 2mu n \bigr\rfloor
\mkern -1mu + \mkern -1mu
\Big\lfloor 2\lfloor\alpha\rfloor \lmyab \alpha \rmyab\mkern 2mu n
+ \lmyab \alpha \rmyab^2\mkern 1mu n \Big\rfloor
\mkern -1mu-\mkern -1mu \Big\lfloor \lfloor\alpha\rfloor \lmyab \alpha \rmyab\mkern 2mu n
+  \lmyab \alpha \rmyab \bigl\lfloor \lmyab \alpha \rmyab\mkern 2mu n \bigr\rfloor\Big\rfloor.
\end{equation}
If we further write
$\lmyab \alpha \rmyab\mkern 2mu n
= \lfloor \lmyab \alpha \rmyab\mkern 2mu n\rfloor
+ \langleat[1.6ex]{70}\mkern-1mu \lmyab \alpha \rmyab\mkern 2mu n\mkern-1mu \rangleat[1.6ex]{70}$, then we can re-write \eqref{new-f} as
\begin{equation*} 
f_\alpha (n)
= \Big\lfloor 2\lfloor\alpha\rfloor
\langleat[1.6ex]{70}\mkern-1mu \lmyab \alpha \rmyab\mkern 2mu n\mkern-1mu \rangleat[1.6ex]{70}
+ \lmyab \alpha \rmyab^2 n\Big\rfloor
- \Big\lfloor \lfloor\alpha\rfloor
\langleat[1.6ex]{70}\mkern-1mu \lmyab \alpha \rmyab\mkern 2mu n\mkern-1mu \rangleat[1.6ex]{70}
+ \lmyab \alpha \rmyab \big\lfloor \lmyab \alpha \rmyab\mkern 2mu n\big\rfloor
\Big\rfloor.
\end{equation*}

\section{Rational parameters $\alpha$}

Recall that if $\alpha\in\mathbb N$, then $f_\alpha(n)=0$ for all $n\in\mathbb N$. Thus, in what follows we shall concentrate on positive rational parameters $\alpha$ which are non-integers. The next lemma states that finding the range of $f_\alpha$ only requires the knowledge of a finite amount of values.

\begin{lemma} \label{cyclic}
Let $a\in\mathbb{N}$ and $b\in\mathbb{N}\mkern-2mu\setminus\mkern-2mu\{1\}$ be relatively prime, and set $\alpha=a/b$. Then
\[
\Range(f_\alpha)
= \{0\}\cup f_\alpha\bigl( \bigl\llbracket 1 , b^2-1 \bigr\rrbracket \bigr).
\]
\end{lemma}

\begin{proof}
Given $n\in\mathbb N$, by the division algorithm, there exist $q, r\in\mathbb Z$ such that
\[
n=b^2q+r,  \quad  r\in \bigl\llbracket 0 , b^2 - 1 \bigr\rrbracket.
\]
Now, if $r=0$, it is easy to see that $f_\alpha (b^2q)=0$. If $r\neq 0$, then
\begin{align*}
f_\alpha(n)
& = \bigg\lfloor a^2q+\frac{ra^2}{b^2}\bigg\rfloor
- \bigg\lfloor\frac{a}{b} \left\lfloor abq + \frac{ar}{b}\right\rfloor \bigg\rfloor\\
& = a^2q + \bigg\lfloor \frac{ra^2}{b^2}\bigg\rfloor - a^2q
- \bigg\lfloor \frac{a}{b} \left\lfloor \frac{ar}{b}\right\rfloor \bigg\rfloor \\
&= f_\alpha(r). \qedhere
\end{align*}
\end{proof}

\subsection{Sub-unitary rational parameters $\alpha$}

We will be assuming here that $\alpha\in\mathbb Q\cap (0, 1)$. The function introduced in the following proposition is reminiscent of the concept of mod operator \cite{K}.
\begin{proposition}
\label{p-s}
Let $a, b \in \mathbb{N}$ be relatively prime with $a<b$. Define the function
\[
\phi_{a,b}:\mathbb{N} \to \mathbb{N}\cup \{0\}
\]
by
\[
\phi_{a,b}\mkern-1mu(n) = n a - \left\lfloor n\frac{a}{b} \right\rfloor b, \quad n \in \mathbb{N}.
\]
Then, $\Range(\phi_{a, b}) = \bigl\llbracket 0 , b-1 \bigr\rrbracket.$

\end{proposition}
\begin{proof}
Using \eqref{defnflce}, we have
\[
na-b<\left\lfloor n\frac{a}{b} \right\rfloor b\leq na,
\]
proving that $\phi_{a,b}\mkern-1mu(n)\in \bigl\llbracket 0 , b-1 \bigr\rrbracket$ for all $n \in \mathbb{N}$. We prove next that $\bigl\llbracket 0 , b-1 \bigr\rrbracket \subseteq \Range(\phi_{a, b})$.
Since $a$ and $b$ are relatively prime, there exist $x,y \in \mathbb{N}$ that satisfy B\'ezout's identity
\[
a x - b y = 1.
\]
Let $r\in \bigl\llbracket 1 , b-1 \bigr\rrbracket$ be arbitrary. Then
\[
\frac{a}{b} xr = yr +  \frac{r}{b}.
\]
Since $0< r/b <1$ and $y \in \mathbb{N}$ we have
\[
\left\lfloor xr \frac{a}{b} \right\rfloor = yr
\]
and
\[
\phi_{a,b}\mkern-1mu(xr) = xra - yrb = r(xa-yb)= r. \qedhere
\]
\end{proof}
The following simple lemma will be needed in our next result.

\begin{lemma}
\label{p-f}
Let $a, b \in \mathbb{N}$ be relatively prime with $a<b$. The range of the function
\begin{equation}\label{eq-fun}
(x,y) \mapsto a x + b y, \qquad x,y\in \bigl\llbracket 0 , b-1 \bigr\rrbracket,
\end{equation}
includes $\bigl\llbracket 0 , b^2 - 1 \bigr\rrbracket$ if and only if $a = 1$.
\end{lemma}

\begin{proof}
Let $a, b \in \mathbb{N}$. Notice that the minimum of the function defined in \eqref{eq-fun} is $0$, its maximum is $(a+b)(b-1)$, so that its codomain is $\bigl\llbracket 0 , (a+b)(b-1)  \bigr\rrbracket$. Note also that the domain has $b^2$ elements while the  codomain has $1+(a+b)(b-1)$ elements, and $b^2 \leq 1+(a+b)(b-1)$. Since
\[
ax_1 + b y_1 = ax_2 + b y_2 \quad \text{with} \quad x_1,y_1, x_2,y_2 \in \bigl\llbracket 0 , b-1 \bigr\rrbracket,
\]
implies the equality of the pairs $(x_1,y_1) = (x_2,y_2)$, the function in \eqref{eq-fun} is an injection. If $a=1$, then the function in \eqref{eq-fun} is a bijection as its domain and its codomain have the same number of elements. Let $a > 1$. Then, by $a<b$, we also get $b>1$. Since the set $\bigl\llbracket 0 , b^2-1 \bigr\rrbracket$ and the range of the function \eqref{eq-fun} have the same number of elements (namely, $b^2$), and since the maximum of the function is
\[
(a+b)(b-1) = b-a + (a-2) b + b^2
\]
which is strictly greater than $b^2 - 1$, it is not possible that $\bigl\llbracket 0 , b^2-1 \bigr\rrbracket$ is a subset of the range. In particular, then $1$ is not in the range of the function defined in \eqref{eq-fun}. This proves the claim.
\end{proof}

\begin{proposition} \label{subunitary}
Let \(\alpha \in \mathbb{Q}\cap (0,1]\). The following dichotomy holds:

\begin{enumerate}
\renewcommand*\theenumi{\roman{enumi}}
\renewcommand*\labelenumi{\rm{(\theenumi)}}

\item \label{subunitary-i1}
If \(\alpha = 1/b\) with \(b \in\mathbb{N}\), then $\Range(f_\alpha)=\{0\}$.

\item \label{subunitary-i2}
If \(\alpha = a/b\) with \(a, b \in\mathbb{N}\) relatively prime and \(1<a < b\), then \(\Range(f_\alpha) = \{0,1\}\).

\end{enumerate}
\end{proposition}
\begin{proof}
By Proposition \ref{largest-range}, we know that $\Range (f_\alpha)\subseteq \{0,1\}$. Moreover, either using Corollary~\ref{01range} or Lemma~\ref{cyclic} (or, further noting that $f_\alpha (kb)=0$ for all $k\in\mathbb N$), we have that \(\{0\}\subseteq \Range (f_\alpha)\).

Now let us assume that $f_\alpha(n) = 0$ for all $n \in \mathbb{N}$. Then, for all $n \in \mathbb{N}$, we have
\begin{equation*} 
\left\lfloor n\frac{a^2}{b^2} \right\rfloor =
\biggl\lfloor \Bigl\lfloor n\frac{a}{b} \Bigr\rfloor  \frac{a}{b} \biggr\rfloor.
\end{equation*}

Let $r \in \bigl\llbracket 0 , b^2-1 \bigr\rrbracket$ be arbitrary. Let $a, b\in\mathbb N$ be relatively prime and such that $a<b$. By Proposition~\ref{p-s} applied to $a^2$ and $b^2$, there exists $m \in \mathbb{N}$ such that $\phi_{a^2,b^2}(m) = r$. Then
\begin{equation}\label{eq-s1}
m\frac{a^2}{b^2} = \left\lfloor m\frac{a^2}{b^2} \right\rfloor + \frac{r}{b^2}.
\end{equation}
By Proposition~\ref{p-s}, we have
\begin{equation}\label{eq-s2}
m\frac{a}{b} = \left\lfloor m\frac{a}{b} \right\rfloor + \frac{\phi_{a,b}\mkern-0.75mu(m)}{b}
\end{equation}
and
\begin{equation}\label{eq-s3}
\left\lfloor m\frac{a}{b} \right\rfloor \frac{a}{b} =
\biggl\lfloor\left\lfloor m\frac{a}{b} \right\rfloor \frac{a}{b} \biggr\rfloor + \frac{\phi_{a,b}\mkern-4mu\left(\left\lfloor m\frac{a}{b} \right\rfloor\right)}{b}.
\end{equation}
Multiplying \eqref{eq-s2} by $a/b$ and using \eqref{eq-s3} we get
\begin{equation}\label{eq-s4}
m\frac{a^2}{b^2} =
\biggl\lfloor\left\lfloor m\frac{a}{b} \right\rfloor \frac{a}{b} \biggr\rfloor + \frac{\phi_{a,b}\mkern-4mu\left(\left\lfloor m\frac{a}{b} \right\rfloor\right)}{b} + \frac{\phi_{a,b}\mkern-0.75mu(m)}{b} \frac{a}{b}.
\end{equation}
Now \eqref{eq-s4}, \eqref{eq-s2} and \eqref{eq-s1} imply
\[
r = a\mkern 1mu \phi_{a,b}\mkern-0.75mu(m) + b\mkern 1mu \phi_{a,b}\mkern-4.5mu\left(\left\lfloor m\tfrac{a}{b} \right\rfloor\right).
\]
Since $\phi_{a,b}\mkern-0.75mu(m),\phi_{a,b}\mkern-4mu\left(\left\lfloor m\tfrac{a}{b} \right\rfloor\right) \in \bigl\llbracket 0 , b-1 \bigr\rrbracket$ and since $r\in \bigl\llbracket 0 , b^2-1 \bigr\rrbracket$ was arbitrary, Lemma~\ref{p-f} implies that $a=1$.
\end{proof}

\begin{remark}\label{rem6} 
Note that, using Lemma~\ref{cyclic}, we can further rephrase the dichotomy in Proposition~\ref{subunitary} as follows: For \(\alpha = a/b\) with \(a, b \in\mathbb{N}\) relatively prime and \(a < b\), the equation
\[
\lfloor \alpha^2 n\rfloor
= \bigl\lceil\alpha \lfloor \alpha n \rfloor \bigr\rceil
\]
is solvable for \(n\) in the set $\bigl\llbracket 0 , b^2-1 \bigr\rrbracket\setminus b \mkern1mu\bigl\llbracket 0 , b-1 \bigr\rrbracket$.
\end{remark}

\begin{remark}
\label{rem7} 
The proof of Proposition~\ref{subunitary} does not explicitly indicate how to exhibit the value 1 in the $\Range(f_\alpha)$ in the case (\ref{subunitary-i2}). We show how this can be achieved.

Let \(a, b \in \mathbb{N}\) be relatively prime and such that \(1 < a < b\). Let \(n, k \in \mathbb{N}\) be such that they satisfy B\'ezout's identity
\begin{equation*} 
a^2 n - b^2 k = 1.
\end{equation*}
Since
\begin{equation*} 
k < \frac{a^2}{b^2} n = k + \frac{1}{b^2} < k+1,
\end{equation*}
we deduce
\begin{equation*} 
\left\lfloor \frac{a^2}{b^2} n \right\rfloor = k.
\end{equation*}
Set
\begin{equation*} 
\left\lfloor \frac{a}{b} n \right\rfloor = r.
\end{equation*}
Then, equivalently,
\begin{equation*} 
r <  \frac{a}{b} n < r+1.
\end{equation*}
Hence,
\(
an - br \in \mathbb{N}.
\)
Therefore, since \(a > 1\), we have
\[
0 < a (an-br) - 1 = a^2 n -abr -a^2 n + b^2 k = - abr  + b^2 k = b( bk - ar ).
\]
Thus, \(l = bk - ar \in \mathbb{N}\), in particular \(ar < bk\), and consequently,
\[
\frac{a}{b} \left\lfloor \frac{a}{b} n \right\rfloor = \frac{a}{b} r = \frac{ar}{b} < \frac{bk}{b} = k.
\]
The last strict inequality yields
\[
\biggl\lfloor \frac{a}{b} \Bigl\lfloor \frac{a}{b} n \Bigr\rfloor \biggr\rfloor < k.
\]
In conclusion,
\[
f_{a/b}(n) = \left\lfloor \frac{a^2}{b^2} n \right\rfloor - \biggl\lfloor \frac{a}{b} \Bigl\lfloor \frac{a}{b} n \Bigr\rfloor \biggr\rfloor>0,
\]
hence $f_{a/b}(n)=1$.
\end{remark}

\subsection{Supra-unitary rational parameters $\alpha$}

Throughout this subsection, we let $\alpha > 1$ be a rational number which is not an integer. Our first result settles the exact range of $f_\alpha$ for supra-unitary rational numbers with denominator 2.

\begin{proposition}
\label{r2}
Let $s \in \mathbb{N}$. Then
\begin{enumerate}
\renewcommand*\theenumi{\roman{enumi}}
\renewcommand*\labelenumi{\rm{(\theenumi)}}

\item \label{r2-i1}
$\Range\bigl(f_{2s+\frac{1}{2}} \bigr) = \{0,s\} = s  \bigl\llbracket 0 , 1 \bigr\rrbracket;$

\item \label{r2-i2}
$\Range\bigl( f_{2s-1+\frac{1}{2}} \bigr) =\{0,s-1,s\} = s \bigl\llbracket 0 , 1 \bigr\rrbracket \cup (s-1) \bigl\llbracket 1 ,1 \bigr\rrbracket.$

\end{enumerate}
\end{proposition}

\begin{proof}
We consider first the case where $\alpha = 2s+\frac{1}{2}$. Then, for all $n\in\mathbb{N}$, we have
\[
\lfloor\alpha n\rfloor=2sn+\Big\lfloor \frac{n}{2}\Big\rfloor = 
\begin{cases}
2sn+\dfrac{n}{2} & \text{if} \quad n \quad \text{is even},    \\[7pt]
2sn+\dfrac{n-1}{2} & \text{if} \quad n \quad \text{is odd}.
\end{cases}
\]
Therefore,
\[
\bigl\lfloor\alpha \lfloor\alpha n\rfloor\bigr\rfloor 
= \begin{cases}
4s^2n+2sn+\Big\lfloor \dfrac{n}{4}\Big\rfloor & \text{if} \quad n \quad \text{is even},  \\[7pt]
4s^2n+2sn-s+\Big\lfloor \dfrac{n-1}{4}\Big\rfloor  & \text{if} \quad n \quad \text{is odd}.
\end{cases} 
\]
Also, $\lfloor\alpha^2 n\rfloor=4s^2n+2sn+\Big\lfloor \frac{n}{4}\Big\rfloor.$ Thus,
\[
\lfloor\alpha^2 n\rfloor - \bigl\lfloor\alpha\lfloor\alpha n\rfloor\bigr\rfloor
=\begin{cases}
0 & \text{if} \quad n \quad \text{is even},  \\[7pt]
s+\Big\lfloor \dfrac{n}{4}\Big\rfloor-\Big\lfloor \dfrac{n-1}{4}\Big\rfloor & \text{if} \quad n \quad \text{is odd}.
\end{cases}
\]
Now, since $\Big\lfloor \dfrac{n}{4}\Big\rfloor-\Big\lfloor \dfrac{n-1}{4}\Big\rfloor=0$ for $n$ odd, we obtain that
\[
\Range\bigl(f_{2s+\frac{1}{2}}\bigr)=\{0, s\}.
\]
Suppose next that $\alpha=2s-1+\frac{1}{2}$. Similar calculations to the ones for the even case above show that, for all $n\in\mathbb N$, we have
\[
\lfloor\alpha n\rfloor = 
2sn+\Big\lfloor\mkern-5mu  -\frac{n}{2}\Big\rfloor=2sn-\Big\lceil \frac{n}{2}\Big\rceil = \begin{cases}
2sn-\dfrac{n}{2} & \text{if} \quad n \quad \text{is even},   \\[7pt]
2sn-\dfrac{n+1}{2} & \text{if} \quad n \quad \text{is odd},
\end{cases}
\]
and
\[
\bigl\lfloor\alpha\lfloor\alpha n\rfloor\bigr\rfloor 
= \begin{cases}
4s^2n-2sn+\Big\lfloor \dfrac{n}{4}\Big\rfloor  & \text{if} \quad n \quad \text{is even}, \\[7pt]
4s^2n-2sn-s+\Big\lfloor \dfrac{n+1}{4}\Big\rfloor  & \text{if} \quad n \quad \text{is odd},
\end{cases}
\]
and
$\lfloor\alpha^2 n\rfloor=4s^2n-2sn+\Big\lfloor \dfrac{n}{4}\Big\rfloor$. Consequently,
\[
\lfloor\alpha^2 n\rfloor - 
\bigl\lfloor\alpha\lfloor\alpha n\rfloor\bigr\rfloor
=\begin{cases}
0  & \text{if} \quad n \quad \text{is even},   \\[7pt]
s+\Big\lfloor \dfrac{n}{4}\Big\rfloor-\Big\lfloor \dfrac{n+1}{4}\Big\rfloor  & \text{if} \quad n \quad \text{is odd}.
\end{cases}
\]
Note now that $\Big\lfloor \dfrac{n}{4}\Big\rfloor-\Big\lfloor \dfrac{n+1}{4}\Big\rfloor=-1$ for $n$ odd such that $n\equiv 3$ (mod 4), while $\Big\lfloor \dfrac{n}{4}\Big\rfloor-\Big\lfloor \dfrac{n+1}{4}\Big\rfloor=0$ for $n\not\equiv 3$ (mod 4). This implies that
\[
\Range\bigl(f_{2s-1+\frac{1}{2}}\bigr)=\{0, s-1, s\}. \qedhere
\]
\end{proof}

In what follows, we explore the general situation in which the (non-integer) rational number $\alpha>1$ has denominator $b\in\mathbb N\setminus\{1\}.$  By the division algorithm, we can express our parameter as
\[
\alpha=(sb+u)+\frac{a}{b},
\]
where $s\in \mathbb{N}\cup \{0\}$, $u\in \bigl\llbracket 0 , b-1 \bigr\rrbracket$ and $a\in \bigl\llbracket 1 , b-1 \bigr\rrbracket$. We have the following decomposition of $f_\alpha$.

\begin{proposition}
\label{decomposition-f}
Let $b\in\mathbb N\setminus\{1\}$, $s\in \mathbb{N}\cup \{0\}$, $u\in \bigl\llbracket 0 , b-1 \bigr\rrbracket$ and $a\in \bigl\llbracket 1 , b-1 \bigr\rrbracket$. Then
\[
f_{sb+u+\frac{a}{b}}=s\mkern 1mu \phi_{a, b}+f_{u+\frac{a}{b}},
\]
where $\phi_{a, b}$ is the function defined in Proposition~{\rm\ref{p-s}}.
\end{proposition}

\begin{proof}
Straightforward calculations give that, for $\alpha=(sb+u)+\frac{a}{b},$ we have
\begin{align*}
\lfloor\alpha^2 n\rfloor&=(sb+u)^2n+2asn+\Big\lfloor \frac{2au}{b}n+\frac{a^2}{b^2}n\Big\rfloor\\
& = (sb+u)^2n+2asn-u^2 n + \Big\lfloor\mkern-3mu \left(\frac{a}{b}+u\right)^2 n\mkern-1mu \Big\rfloor.
\end{align*}
Next, we have
\begin{align*}
\alpha\lfloor\alpha n\rfloor&=(sb+u)^2n+asn+sb\Big\lfloor \frac{a}{b}n\Big\rfloor+\frac{au}{b}n+\left(\frac{a}{b}+u\right)\Big\lfloor \frac{a}{b}n\Big\rfloor\\
&=(sb+u)^2n+asn+sb\Big\lfloor \frac{a}{b}n\Big\rfloor+\frac{2au}{b}n+\frac{a^2}{b^2}n-\left(\frac{a}{b}+u\right)
 \langleat[2.5ex]{70} \frac{a}{b}n\mkern-1mu \rangleat[2.5ex]{70}\\
&=(sb+u)^2n+asn+sb\Big\lfloor \frac{a}{b}n\Big\rfloor-u^2 n \\
& \mkern 140mu +\left(\frac{a}{b}+u\right)^2 n  -
\left(\frac{a}{b}+u\right) 
\langleat[2.5ex]{70}\mkern-1mu \left(\frac{a}{b}+u\right)n\mkern-1mu \rangleat[2.5ex]{70}\\
& = (sb+u)^2n+asn+sb\Big\lfloor \frac{a}{b}n\Big\rfloor-u^2 n + 
\left(\frac{a}{b}+u\right) \Big\lfloor\mkern-3mu \left(\frac{a}{b}+u\right) n \mkern-1mu\Big\rfloor.
\end{align*}
Thus,
\[
\bigl\lfloor\alpha\lfloor\alpha n\rfloor\bigr\rfloor=(sb+u)^2n+asn+sb\Big\lfloor \frac{a}{b}n\Big\rfloor
- u^2 n + 
\bigg\lfloor \mkern-5mu
\left(\frac{a}{b} + u\right) 
\Big\lfloor \left(\frac{a}{b}+u\right) n \Big\rfloor\bigg\rfloor
\]
and
\begin{align*}
f_\alpha (n) 
& = asn
- sb\Big\lfloor \frac{a}{b}n\Big\rfloor+\Big\lfloor \left(\frac{a}{b}+u\right)^2 n\Big\rfloor 
- \bigg\lfloor\mkern-3mu\left(\frac{a}{b}+u\right)\Big\lfloor \left(\frac{a}{b}+u\right) n\Big\rfloor\bigg\rfloor\\
& = s\mkern 1mu \phi_{a, b}(n)+f_{u+\frac{a}{b}}(n).\qedhere
\end{align*}

\end{proof}
Two immediate consequences of the decomposition in Proposition~\ref{decomposition-f} are the following.

\begin{corollary}
\label{general-inclusion}
Let $b\in\mathbb N\setminus\{1\}$, $s\in \mathbb{N}\cup \{0\}$, $u\in \bigl\llbracket 0 , b-1 \bigr\rrbracket$ and $a\in \bigl\llbracket 1 , b-1 \bigr\rrbracket$ relatively prime to $b$. Then
\[
\Range \bigl(f_{sb+u+\frac{a}{b}}\bigr)\subseteq s \bigl\llbracket 0 , b-1 \bigr\rrbracket + \bigl\llbracket 0 , u+1 \bigr\rrbracket.
\]
\end{corollary}

\begin{proof}
Noting that $\Big\lceil u+\frac{a}{b}\Big\rceil=u+1$, the inclusion follows from Proposition~\ref{p-s} and Proposition~\ref{largest-range}
\end{proof}

\begin{corollary}
\label{rb0}
Let $b\in\mathbb{N}\mkern-2mu\setminus\mkern-2mu\{1\}$ and $s\in\mathbb N$. Then
\[
\Range\bigl( f_{sb+\frac{1}{b}} \bigr)=s \bigl\llbracket 0 , b-1 \bigr\rrbracket.
\]
\end{corollary}

\begin{proof}
By Proposition~\ref{decomposition-f} and Proposition~\ref{subunitary}, part (i), for all $n\in\mathbb N$ we have
\[
f_{sb+\frac{1}{b}}(n)=s\phi_{1, b}(n)+f_{\frac{1}{b}}(n)=s\phi_{1, b}(n).
\]
The equality \(\Range\bigl( f_{sb+\frac{1}{b}} \bigr) = s \bigl\llbracket 0 , b-1 \bigr\rrbracket\) now follows from Proposition~\ref{p-s}.
\end{proof}

While Proposition~\ref{decomposition-f} provides a good starting point for identifying the exact range of functions such as $f_{sb+\frac{1}{b}}$, we need a slightly modified version of it to deal with other supra-unitary parameters of the form $sb+u+\frac{1}{b}$ with $u>0$.

\begin{lemma} \label{fbu}
Let $b\in\mathbb N\setminus\{1\}$, $s\in \mathbb{N}\cup \{0\}$, and $u\in \bigl\llbracket 0 , b-1 \bigr\rrbracket$. Let $n=bt+r$ with $t, r\in \bigl\llbracket 0 , b-1 \bigr\rrbracket$. Then
\[
f_{sb+u+\frac{1}{b}}(n)=sr+\Big\lfloor \frac{2ur}{b}+\frac{t}{b}+\frac{r}{b^2}\Big\rfloor-\Big\lfloor \frac{ur}{b}+\frac{t}{b}\Big\rfloor. \]
\end{lemma}
In particular, for $u=0$, we obtain
\(
f_{sb+\frac{1}{b}} (n)=sr,
\)
which implies the statement of Corollary~\ref{rb0}.

\begin{proof}
The calculations are similar to the ones in the proof of Proposition~\ref{decomposition-f}. We only present the essential identities. With $\alpha=bs+u+\frac{1}{b}$, we have
\[
\lfloor\alpha^2 n\rfloor=(bs+u)^2 n+2sn+2ut+\Big\lfloor \frac{2ur}{b}+\frac{t}{b}+\frac{r}{b^2}\Big\rfloor
\]
and
\[
\bigl\lfloor\alpha\lfloor\alpha n\rfloor \bigr\rfloor = (bs+u)^2 n+sn+bst+2ut+\Big\lfloor \frac{ur}{b}+\frac{t}{b}\Big\rfloor.
\]
Therefore,
\begin{align*}
f_\alpha(n) &= s(n-bt)
+ \Big\lfloor\frac{2ur}{b}+\frac{t}{b}+\frac{r}{b^2}\Big\rfloor - \Big\lfloor \frac{ur}{b}+\frac{t}{b}\Big\rfloor\\
& = sr + \Big\lfloor\frac{2ur}{b}+\frac{t}{b}+\frac{r}{b^2}\Big\rfloor - \Big\lfloor \frac{ur}{b}+\frac{t}{b}\Big\rfloor.
 \qedhere
\end{align*}
\end{proof}

Our next result states an inclusion relation that improves the inclusion stated in Corollary~\ref{general-inclusion} when $a=1$.

\begin{proposition}
\label{rbu-inclusion}
For all \(b\in \mathbb{N}\mkern-2mu\setminus\mkern-2mu\{1\}\), all $s\in \mathbb{N}\cup \{0\}$, and all $u\in \bigl\llbracket 0 , b-1 \bigr\rrbracket$ such that \(sb+u \in \mathbb{N}\), we have
\[
\Range\bigl( f_{sb+u+\frac{1}{b}} \bigr) \subseteq s \bigl\llbracket 0 , b-1 \bigr\rrbracket + \bigl\llbracket 0 , u \bigr\rrbracket.
\]
\end{proposition}

\begin{proof}
Let $\alpha=sb+u+\frac{1}{b}$. Using the notation from Lemma~\ref{fbu}, with $n=bt+r$ and $t, r\in \bigl\llbracket 0 , b-1 \bigr\rrbracket,$ we can rewrite
\begin{equation}
\label{nice-alternate}
f_{\alpha}(n)=sr+\Big\lfloor \frac{2ur}{b}+\frac{n}{b^2}\Big\rfloor-\Big\lfloor \frac{(ub-1)r}{b^2}+\frac{n}{b^2}\Big\rfloor.
\end{equation}
Let us denote for simplicity
\[
x(n, u)=\frac{2ur}{b}+\frac{n}{b^2}\,\,\text{and}\,\,y(n, u)=\frac{(ub-1)r}{b^2}+\frac{n}{b^2}.
\]
Clearly, if $u=0$, we have $\lfloor x(n, u)\rfloor=\lfloor y(n,u)\rfloor=0.$ If $u\geq 1$, we see that
\[
x(n, u)-y(n, u)=\frac{(ub+1)r}{b^2}\geq 0
\]
and
\[
x(n, u)-y(n, u)\leq \frac{(ub+1)(b-1)}{b^2}=u-\frac{u-1}{b}-\frac{1}{b^2}<u.
\]
Therefore,
\[
0 \leq \bigl\lfloor x(n, u)-y(n,u) \bigr\rfloor
\leq \bigl\lfloor x(n, u) \bigr\rfloor - \bigl\lfloor y(n,u) \bigr\rfloor,
\]
and
\[
\bigl\lfloor x(n, u)\bigr\rfloor - \bigl\lfloor y(n,u) \bigr\rfloor \leq \lfloor x(n, u)-y(n,u)\rfloor+1\leq (u-1)+1=u.
\]
In other words, for all $n\in \bigl\llbracket 0 , b^2-1 \bigr\rrbracket$ and $u\in \bigl\llbracket 0 , b-1 \bigr\rrbracket$, we have
\[
\bigl\lfloor x(n, u) \bigr\rfloor-\bigl\lfloor y(n,u)\bigr\rfloor\in \bigl\llbracket 0 , u \bigr\rrbracket.
\]
This finishes the proof.
\end{proof}
In particular, when $u=1$, we get the following.
\begin{corollary}
\label{u1}
For all \(b\in\mathbb{N}\mkern-2mu\setminus\mkern-2mu\{1\}\) and all $s\in \mathbb{N}\cup \{0\}$ we have
\[
\Range\bigl(f_{sb+1+\frac{1}{b}}\bigr)\subseteq s \bigl\llbracket 0 , b-1 \bigr\rrbracket + \bigl\llbracket 0 , 1 \bigr\rrbracket.
\]
\end{corollary}

\begin{corollary}
\label{s0u1}
Let $b\in\mathbb{N}\mkern-2mu\setminus\mkern-2mu\{1\}$. Then
\[
\Range\bigl(f_{1+\frac{1}{b}}\bigr)=\{0, 1\}.
\]
\end{corollary}
\begin{proof}
Using Lemma~\ref{cyclic} and Corollary~\ref{u1}, we only need to prove that
$1\in\Range\bigl( f_{1+\frac{1}{b}} \bigr).$ We claim that
\(
f_{1+\frac{1}{b}}(b-1)=1.
\)
Indeed, letting $s=0$, $u=1$, $t=0$ and $r=b-1$ in Lemma~\ref{fbu} (or, alternately, in the formula \eqref{nice-alternate}), we have
\[
f_{1+\frac{1}{b}}(b-1)=\Big\lfloor 2-\frac{b+1}{b^2}\Big\rfloor-\Big\lfloor \frac{b-1}{b}\Big\rfloor.
\]
But since $1<2-\frac{b+1}{b^2}<2$ and $0<\frac{b-1}{b}<1$, we get $f_{1+\frac{1}{b}}(b-1)=1.$
\end{proof}

A refinement of the argument in Corollary~\ref{s0u1} gives the next result.
\begin{corollary}
\label{s1u1}
Let $b\in\mathbb{N}\mkern-2mu\setminus\mkern-2mu\{1\}$. Then
\[
\Range\bigl( f_{b+1+\frac{1}{b}} \bigr) = \bigl\llbracket 0 , b \bigr\rrbracket.
\]
\end{corollary}
\begin{proof}
Using again Lemma~\ref{cyclic} and Corollary~\ref{u1}, it suffices to show that $\bigl\llbracket 1 , b \bigr\rrbracket \subseteq\Range\bigl( f_{b+1+\frac{1}{b}} \bigr).$ We proceed in two steps. First, we claim that $f_{b+1+\frac{1}{b}}(b^2-b+1)=1$. Indeed, letting $s=1$, $u=1$, $t=b-1$ and $r=1$ in Lemma~\ref{fbu}, and since $1<1+\frac{b+1}{b^2}<2$, we get
\[
f_{b+1+\frac{1}{b}}\bigl(b^2-b+1\bigr) = 1+\Big\lfloor 1+\frac{b+1}{b^2}\Big\rfloor-1=1.
\]
Secondly, let us fix $r\in \bigl\llbracket 1 , b-1 \bigr\rrbracket$. We will show that there exists a $t(r)\in \bigl\llbracket 0 , b-1 \bigr\rrbracket$ such that $f_{b+1+\frac{1}{b}}\bigl(bt(r)+r\bigr)=r+1$. Indeed, let $s=u=1$ and $t=t(r):=b-1-r$ in Lemma~\ref{fbu} and notice that
\[
f_{b+1+\frac{1}{b}}\bigl(bt+r\bigr) = r
+ \Big\lfloor\mkern-1mu \frac{r}{b}+\frac{r}{b^2} + \frac{b-1}{b}\mkern-1mu\Big\rfloor - \Big\lfloor\mkern-1mu \frac{b-1}{b} \mkern-1mu\Big\rfloor
=
r+\Big\lfloor\mkern-1mu 1+\frac{r-1}{b}+\frac{r}{b^2}\mkern-1mu\Big\rfloor.
\]
Next, observe that \(0<\frac{r-1}{b}+\frac{r}{b^2}<1\). The left-side inequality is trivial while the right-side inequality is equivalent to \((r-1)b+r<b^2\Leftrightarrow r<b\). This proves the claim $f_{b+1+\frac{1}{b}}\bigl(b(b-1-r)+r\bigr)=r+1$ and concludes our proof.
\end{proof}

The following corollary generalizes Corollaries~\ref{s0u1} and~\ref{s1u1}.

\begin{corollary} \label{ssu1}
For all $b\in\mathbb N\mkern-2mu\setminus\mkern-2mu\{1\}$ and $s\in \mathbb{N}\cup \{0\}$ we have
\[
\Range\bigl(f_{sb+1+\frac{1}{b}}\bigr)
= \bigl( s \bigl\llbracket 0 , b-1 \bigr\rrbracket \bigr) \cup
\bigl(s \bigl\llbracket 1 , b-1 \bigr\rrbracket + \{1\}\bigr).
\]
\end{corollary}
\begin{proof}
It is sufficient to consider \(s> 1\). By Lemma~\ref{cyclic} and Corollary~\ref{u1}, it suffices to show that
$1\not\in\Range\bigl( f_{sb+1+\frac{1}{b}} \bigr)$ and $s \bigl\llbracket 1 , b-1 \bigr\rrbracket + \{0,1\} \subseteq\Range\bigl(f_{sb+1+\frac{1}{b}}\bigr)$. First of all, using the notation from Lemma~\ref{fbu}, for all $r\in \bigl\llbracket 1 , b-1 \bigr\rrbracket$, we have $f_{sb+1+\frac{1}{b}}(bt+r)\geq sr\geq s > 1$ (and, as we have seen before, $f_{sb+1+\frac{1}{b}}(bt)=0$). For the second claim, we will show that for \emph{fixed}
$r\in \bigl\llbracket 1, b-1 \bigr\rrbracket$, we can find $t_0(r), t_1(r)\in \bigl\llbracket 0 , b-1 \bigr\rrbracket$ such that
\[
f_{\alpha}\bigl(bt_0(r)+r \bigr) = sr\quad \text{and}\quad
f_{\alpha}\bigl( bt_1(r)+r \bigr) = sr+1.
\]
Equivalently, we want to show that, for fixed $r\in \bigl\llbracket 1 , b-1 \bigr\rrbracket$, the equations
\begin{equation}
\label{rt0}
\Big\lfloor \frac{2r}{b}+\frac{r}{b^2}+\frac{t}{b}\Big\rfloor=\Big\lfloor \frac{r}{b}+\frac{t}{b}\Big\rfloor,
\end{equation}
respectively
\begin{equation}
\label{rt1}
\Big\lfloor \frac{2r}{b}+\frac{r}{b^2}+\frac{t}{b}\Big\rfloor=\Big\lceil \frac{r}{b}+\frac{t}{b}\Big\rceil
\end{equation}
are solvable in $t\in \bigl\llbracket 0 , b-1 \bigr\rrbracket$.

The equation ~\eqref{rt0} has the solution $t_0=t_0(r)=b-r$. Indeed, $\displaystyle\Big\lfloor \frac{r}{b}+\frac{t_0}{b}\Big\rfloor=1$ and
\[
\frac{2r}{b}+\frac{r}{b^2}+\frac{t_0}{b}=\frac{r}{b}+\frac{r}{b^2}+1\in (1, 2)
\]
since $r\leq b-1<\displaystyle\frac{b^2}{b+1}$, therefore also $\displaystyle\Big\lfloor \frac{2r}{b}+\frac{r}{b^2}+\frac{t_0}{b}\Big\rfloor=1$. The equation ~\eqref{rt1} has the solution $t_1(r)=b-1-r$ by the argument proving Corollary~\ref{s1u1}.
\end{proof}
Similar considerations to the ones proving Proposition~\ref{rbu-inclusion} and Corollary~\ref{ssu1} give the following result.

\begin{corollary}
\label{ssu-1}
Let $b\in\mathbb{N}\mkern-3mu\setminus\mkern-3mu\{1\}$ and $s\in\mathbb N$. Then
\[
\Range\bigl(f_{sb-1 +\frac{1}{b}}\bigr)
= \bigl(s \bigl\llbracket 0 , b-1 \bigr\rrbracket \bigr) \cup
\bigl(s \bigl\llbracket 1 , b-1 \bigr\rrbracket - \{1\}\bigr).
\]
\end{corollary}

\begin{proof}
Let $\alpha=sb-1+\frac{1}{b}=(s-1)b+(b-1)+\frac{1}{b}$. With $s$ replaced by $s-1\in \mathbb{N}\cup \{0\}$, and $u=b-1$, in the proof of Proposition~\ref{rbu-inclusion}, we find that \eqref{nice-alternate} can be rewritten as
\[
f_\alpha (n)=(s-1)r+\Big\lfloor \frac{2(b-1)r}{b}+\frac{n}{b^2}\Big\rfloor-\Big\lfloor \frac{((b-1)b-1)r}{b^2}+\frac{n}{b^2}\Big\rfloor.
\]
As before, letting now $n=bt+r$ with $t, r\in \bigl\llbracket 0 , b-1 \bigr\rrbracket,$ a few simple manipulations yield
\begin{equation}
\label{nice-alternate2}
f_\alpha (bt+r)=sr -\Bigl(\Big\lfloor\frac{r}{b}+\frac{t}{b}\Big\rfloor-\Big\lfloor -\frac{2r}{b}+\frac{r}{b^2}+\frac{t}{b}\Big\rfloor\Bigr).
\end{equation}
Denote
\[
x(b, t, r)=\frac{r}{b}+\frac{t}{b}\,\,\text{and}\,\,y(b, t, r)=-\frac{2r}{b}+\frac{r}{b^2}+\frac{t}{b}.
\]
We have
\[
x(b, t, r)-y(b, t, r)=\frac{(b-1)r}{b^2}\in \Bigl(0, \frac{(b-1)^2}{b^2}\Bigr]\subseteq (0, 1).
\]
This implies that $\bigl\lfloor x(b, t, r)-y(b, t, r)\bigr\rfloor =0$, thus
\begin{equation}
\label{nice3}
\bigl\lfloor x(b, t, r)\bigr\rfloor - \bigl\lfloor y(b, t, r)\bigr\rfloor \in \{0, 1\}.
\end{equation}
Using \eqref{nice-alternate2} and \eqref{nice3}, we conclude that
\[
\Range\bigl(f_{sb-1+\frac{1}{b}}\bigr) \subseteq
s \bigl\llbracket 0 , b-1 \bigr\rrbracket - \{0, 1\}.
\]
Since $-1\not\in\Range\bigl( f_{sb-1+\frac{1}{b}}\bigr)$, we have
\[
\Range\bigl( f_{sb-1 +\frac{1}{b}} \bigr) \subseteq \bigl(s \bigl\llbracket 0 , b-1 \bigr\rrbracket \bigr) \cup \bigl(s \bigl\llbracket 1 , b-1 \bigr\rrbracket - \{1\}\bigr).
\]
For the converse inclusion, as in the proof of Corollary \ref{ssu1}, it suffices to show that
\[
s \bigl\llbracket 1 , b-1 \bigr\rrbracket - \{0,1\} \subseteq
\Range\bigl( f_{sb-1+\frac{1}{b}} \bigr),
\]
which is equivalent to showing that for \emph{fixed}
$r\in \bigl\llbracket 1 , b-1 \bigr\rrbracket$, we can find $t(r), \tilde t(r)\in \bigl\llbracket 0 , b-1 \bigr\rrbracket$ such that
\[
f_{\alpha}\bigl(bt(r)+r \bigr) = sr\quad \text{and}\quad
f_{\alpha}\bigl( b\tilde t(r)+r \bigr) = sr-1.
\]
Equivalently, we want to show that, for fixed $r\in \bigl\llbracket 1, b-1 \bigr\rrbracket$, the equations
\begin{equation} \label{rt}
\Big\lfloor\mkern-4mu -\frac{r}{b}+\frac{t}{b} \Big\rfloor 
= \Big\lfloor\mkern-4mu -\frac{2r}{b}+\frac{r}{b^2}+\frac{t}{b} \Big\rfloor,
\end{equation}
respectively
\begin{equation}
\label{rtt}
\Big\lfloor\mkern-4mu -\frac{r}{b}+\frac{t}{b}\Big\rfloor 
= 1+\Big\lfloor\mkern-4mu -\frac{2r}{b}+\frac{r}{b^2}+\frac{t}{b}\Big\rfloor
\end{equation}
are solvable in $t\in \bigl\llbracket 0 , b-1 \bigr\rrbracket$. The equation \eqref{rt} has the solution $t(r)=r-1$. The equation \eqref{rtt} has the solution $\tilde t(r)=r$.
\end{proof}

The exact range results stated in Proposition~\ref{r2}, Corollary~\ref{rb0}, Corollary~\ref{ssu1} and Corollary~\ref{ssu-1} capture \emph{some} generic supra-unitary (non-integer) rational parameters $\alpha$ with denominator $b\geq 2$. Specifically, we understand precisely the cases where $\alpha=k+\frac{1}{b}$ and $k\equiv u$ (mod $b$) with $u\in\{0, 1, b-1\}$, but, even in the small denominator cases, we are missing lots of the rational supra-unitary parameters. The result below is a counterpart to Proposition~\ref{r2} which settles the exact ranges for \emph{all} supra-unitary rational parameters with denominator 3.

\begin{proposition} \label{r3}
Let $s\in\mathbb N$. Then
\begin{enumerate}
\renewcommand*\theenumi{\roman{enumi}}
\renewcommand*\labelenumi{\rm{(\theenumi)}}

\item \label{r3-i1}
$\Range\bigl(f_{3s+\frac{1}{3}}\bigr) = \{0, s, 2s\};$

\item \label{r3-i2}
$\Range\bigl( f_{3s+\frac{2}{3}} \bigr) = \{0, s, 2s, 2s+1\};$

\item \label{r3-i3}
$\Range\bigl( f_{3s-2+\frac{1}{3}} \bigr)=\{0, s-1, s, 2s-2, 2s-1\};$ 

\item \label{r3-i4}
$\Range\bigl( f_{3s-2+\frac{2}{3}} \bigr) = \{0, s-1, s, 2s-1\};$

\item \label{r3-i5}
$\Range\bigl( f_{3s-1+\frac{1}{3}} \bigr) =\{0, s-1, s, 2s-1, 2s\};$

\item \label{r3-i6}
$\Range\bigl( f_{3s-1+\frac{2}{3}} \bigr)=\{0, s-1, s, 2s-1, 2s\}.$
\end{enumerate}
\end{proposition}

\begin{proof}
(\ref{r3-i1}) follows from Corollary~\ref{rb0}, (\ref{r3-i3}) follows from Corollary~\ref{ssu1} (via the substitution $s\to s-1$) and (\ref{r3-i5}) follows from Corollary~\ref{ssu-1}.

(\ref{r3-i2}) Let $\alpha=3s+\frac{2}{3}$. Based on Lemma~\ref{cyclic}, we know that $0\in\Range(f_\alpha)$ and for all the other values in $\Range (f_\alpha)$ we can concentrate on computing $f_\alpha (n)$ with $n\in \bigl\llbracket 0 , 8 \bigr\rrbracket$. For such an $n$, we have
\[
\alpha^2 n=9s^2n+4sn+\frac{4n}{9}\Rightarrow \lfloor\alpha^2 n\rfloor=9s^2n+4sn+\Big\lfloor\frac{4n}{9}\Big\rfloor.
\]
Similarly, we compute
\[
\bigl\lfloor\alpha\lfloor\alpha n\rfloor\bigr\rfloor =
9s^2n+2sn +
3s\Big\lfloor\frac{2n}{3}\Big\rfloor +
\bigg\lfloor\frac{2}{3}\Big\lfloor\frac{2n}{3}
\Big\rfloor\bigg\rfloor.
\]
Hence,
\begin{equation}
\label{caseb}
f_\alpha(n) = 2sn - 3s\Big\lfloor\frac{2n}{3}\Big\rfloor + \Big\lfloor\frac{4n}{9}\Big\rfloor -
\bigg\lfloor\frac{2}{3}\Big\lfloor\frac{2n}{3}
\Big\rfloor\bigg\rfloor.
\end{equation}
Using \eqref{caseb}, we compute:
\begin{multline*}
f_\alpha(3)=f_\alpha(6)=0, \quad f_{\alpha}(2)=f_\alpha(5)=f_\alpha(8)=s, \\ f_\alpha(1)=f_\alpha(4)=2s, \quad f_\alpha(7)=2s+1.
\end{multline*}

(\ref{r3-i4}) Let $\alpha=3s-2+\frac{2}{3}=3s-\frac{4}{3}$. After some algebraic manipulations, we obtain
\begin{equation}
\label{cased}
f_\alpha(n) = - sn + 3s\Big\lceil\frac{n}{3}\Big\rceil + \Big\lfloor\frac{7n}{9}\Big\rfloor -
\Big\lceil\frac{n}{3}\Big\rceil
- \bigg\lfloor\frac{n}{3} + \frac{1}{3}\Big\lceil\frac{n}{3}
\Big\rceil\bigg\rfloor.
\end{equation}
Using \eqref{cased}, we can now compute:
\begin{multline*}
f_\alpha(3)=f_\alpha(6)=0, \quad f_{\alpha}(2)=f_\alpha(5)=s-1, \\
f_\alpha(8)=s, \quad f_\alpha(1)=f_\alpha(4)=f_\alpha(7)=2s-1.
\end{multline*}

(\ref{r3-i6}) Let $\alpha=3s-1+\frac{2}{3}=3s-\frac{1}{3}$. Since $n\in \bigl\llbracket 0 , 8 \bigr\rrbracket$, we have
\begin{equation}
\label{casef}
f_\alpha(n) = -sn + 3s\Big\lceil\frac{n}{3}\Big\rceil -
\bigg\lfloor\frac{1}{3}\Big\lceil\frac{n}{3}
\Big\rceil\bigg\rfloor.
\end{equation}
Using \eqref{casef}, we obtain:
\begin{multline*}
f_\alpha(3)=f_\alpha(6)=0, \quad f_\alpha(8)=s-1, \quad f_\alpha(2)=f_\alpha(5)=s, \\
f_\alpha(7)=2s-1, \quad f_\alpha(1)=f_\alpha(4)=2s,  
\end{multline*}
completing the proof.  
\end{proof}

Besides the general cases already settled above and performing ad-hoc calculations in the small denominators cases, such as in Proposition~\ref{r2} or Proposition ~\ref{r3}, identifying the \emph{exact} $\Range(f_\alpha)$ for supra-unitary rational values of $\alpha$ seems to be a rather difficult task. It is not clear to these authors how to proceed \emph{in general} for such parameters beyond the situations already addressed above. For example, our Wolfram Mathematica experimentations, see \cite{BC1}, are pretty convincing that the following statement should be true, but we are unable to prove it.

\begin{conjecture}
\label{ssu-1a}
For all $b\in\mathbb{N}\mkern-3mu\setminus\mkern-3mu\{1\}$, all \(a \in \bigl\llbracket 1 , b-1 \bigr\rrbracket\) relatively prime to \(b\), and all $s\in\mathbb N$ we have
\[
\Range\bigl( f_{sb - \frac{a}{b}} \bigr)
= \bigl(s \bigl\llbracket 0 , b-1 \bigr\rrbracket \bigr) \cup \bigl(s \bigl\llbracket 1 , b-1 \bigr\rrbracket - \{1\}\bigr).
\]
\end{conjecture}

If indeed true, Conjecture~\ref{ssu-1a} would generalize Corollary~\ref{ssu-1} which corresponds to $a=b-1$. The following lemma is in the spirit of Lemma~\ref{fbu}. The calculations are similar to the ones done there and left to the interested reader. We include it here since it may indicate a possible direction to deal with cases such as those in Conjecture~\ref{ssu-1a}.

\begin{lemma} \label{fbuv}
Let $b\in\mathbb{N}\mkern-1mu\setminus\mkern-2mu\{1\}$, $u\in \bigl\llbracket 0 , b-1 \bigr\rrbracket$ and $a\in \bigl\llbracket 1 , b-1 \bigr\rrbracket$. Let $n=bt+r,$ with $t, r\in \bigl\llbracket 0 , b-1 \bigr\rrbracket$. Then
\[
f_{sb+u+\frac{a}{b}}(n) = sar
- (sb+u)\Big\lfloor\mkern-2mu\frac{ar}{b}\mkern-2mu\Big\rfloor
+ \bigg\lfloor\mkern-2mu \frac{2uar}{b}+\frac{a^2t}{b} + \frac{a^2r}{b^2}\mkern-2mu\bigg\rfloor
- \bigg\lfloor\mkern-2mu \frac{uar}{b}+\frac{a^2t}{b} + \frac{a}{b}\Big\lfloor\frac{ar}{b}\Big\rfloor\mkern-2mu\bigg\rfloor.
\]
\end{lemma}

\section{Irrational parameters $\alpha$}

The goal of this section is to discuss a generalization of \cite{XX} and present several conjectures inspired by it and from extensive experimentations in Wolfram Mathematica, see \cite{BC1}.

\subsection{Generalized golden ratios} Given $t\in\mathbb N$ such that
\[
t\not\in\bigl\{(s-1)s: s\in\mathbb N\bigr\},
\]
we let
\begin{equation}
\label{alpha-special}
\alpha (t) := \frac{1+\sqrt{1+4t}}{2};
\end{equation}
note that $\alpha (t)\in (\mathbb R\mkern-2mu\setminus\mkern-2mu\mathbb Q)\cap (1, \infty)$. The number $\alpha (1)=\varphi$ is the famous golden ratio; $\alpha(504)=\tau$ is the irrational number in \cite{XX}. Our main goal is to prove the following.

\begin{proposition} \label{golden}
For all $t\in\mathbb N$ and $\alpha(t)$ as in \eqref{alpha-special}, we have
$$
\Range(f_{\alpha (t)}) = \bigl\llbracket 1 , \lfloor \alpha (t) \rfloor \bigr\rrbracket .
$$
\end{proposition}

\begin{proof}
Fix $t\in\mathbb N$ and simply write $\alpha (t)=\alpha$ for ease of notation. Let $n\in\mathbb N$. We begin by observing that $\alpha^2-\alpha=t$, hence $\alpha^2 n=tn+\alpha n$. Therefore
\[
f_{\alpha}(n)=tn+\lfloor\alpha n\rfloor-\bigl\lfloor \alpha  \lfloor \alpha n \rfloor \bigr\rfloor.
\]
Our first goal is to show that $0<f_\alpha (n)<\alpha$. Since $f_\alpha (n)\in\mathbb Z$, this in turn will give that $1\leq f_\alpha(n)\leq\lfloor\alpha\rfloor$ for all $n\in\mathbb N$, and hence the inclusion of sets
\[
\Range(f_{\alpha (t)})\subseteq \bigl\llbracket 1 , \lfloor \alpha (t) \rfloor \bigr\rrbracket.
\]
Note that, by \eqref{defnflce}, we have
\begin{align*}
f_\alpha (n)
& < tn+\lfloor\alpha n\rfloor-(\alpha \lfloor\alpha n\rfloor-1) \\
& = tn+1-(\alpha-1)\lfloor\alpha n\rfloor\\
& = tn+1-(\alpha-1)(\alpha n- \lmyab \alpha n \rmyab) \\
& = n(t-\alpha(\alpha-1))+1+(\alpha-1) \mkern 1mu \lmyab \alpha n \rmyab\\
& = 1+(\alpha-1) \mkern 1mu \lmyab \alpha n \rmyab \\
& < 1+(\alpha-1) \\
& =\alpha.
\end{align*}
Similarly, we can write
\begin{align*}
f_\alpha (n)
& > tn+\lfloor\alpha n\rfloor-\alpha \lfloor\alpha n\rfloor \\
& = tn-(\alpha-1)\lfloor\alpha n\rfloor\\
& = tn-(\alpha-1)(\alpha n-\lmyab \alpha n \rmyab) \\
& = (\alpha-1) \mkern 1mu \lmyab \alpha n \rmyab \\
& >0;
\end{align*}
the last strict inequality uses the fact that $\alpha>1$ and $\alpha\not\in\mathbb Q$ (that is, $\lmyab \alpha n \rmyab = 0$ is equivalent to $\alpha=\lfloor\alpha n \rfloor/n\in\mathbb Q$). We record for future purposes that we obtained the stronger inequalities
\begin{equation} \label{ineq}
(\alpha-1)\mkern 1mu \lmyab \alpha n \rmyab <f_\alpha (n)<1+(\alpha-1) \mkern 1mu\lmyab \alpha n \rmyab,\quad \forall n\in\mathbb N.
\end{equation}

We will prove next the inclusion
\[
\Range(f_{\alpha (t)})\supseteq \bigl\llbracket 1 , \lfloor \alpha (t) \rfloor \bigr\rrbracket.
\]
The essential tool in this regard is Kronecker's approximation theorem \cite{W} which states that for irrational $\alpha$, the set $\bigl\{\lmyab \alpha n \rmyab: n\in\mathbb N\bigr\}$ is dense in $[0, 1]$. Let $m \in \bigl\llbracket 1 , \lfloor \alpha \rfloor \bigr\rrbracket.$ Since $m/\alpha\in (0, 1)$, we have that
\begin{equation}
\label{dense}
\forall\mkern 1.5mu \epsilon>0 \quad \exists\,n\in\mathbb N: \quad \frac{m}{\alpha}-\epsilon< \lmyab \alpha n \rmyab < \frac{m}{\alpha}+\epsilon.
\end{equation}
Let $\epsilon = \lmyab \alpha\rmyab/{t}\in (0, 1)$. We will break our discussion into two cases. If $m=\lfloor\alpha\rfloor$, by \eqref{dense}, there exists some $n\in\mathbb N$ such that
\[
(\alpha-1) \mkern 1mu \lmyab \alpha n \rmyab
> (\alpha-1)\frac{\lfloor\alpha\rfloor}{\alpha}
- (\alpha-1)\frac{\lmyab \alpha\rmyab}{\alpha(\alpha-1)}.
\]
By \eqref{ineq}, we get that
\[
f_\alpha (n)
>
(\alpha-1) \frac{\lfloor\alpha\rfloor}{\alpha}
- \frac{\lmyab \alpha\rmyab}{\alpha}
=
\lfloor\alpha\rfloor
- \frac{\lfloor\alpha\rfloor+\lmyab \alpha\rmyab}{\alpha}=\lfloor\alpha\rfloor-1.
\]
Since $f_\alpha (n)\leq \lfloor\alpha\rfloor$, we get that $f_\alpha (n)= \lfloor\alpha\rfloor$.

If $1\leq m\leq \lfloor\alpha\rfloor -1$, since $\epsilon<\frac{1}{\alpha(\alpha-1)}$, from \eqref{dense} we get some $n\in\mathbb N$ such that
\[
\frac{m(\alpha-1)}{\alpha}-\frac{1}{\alpha}
< (\alpha-1)\mkern 1mu \lmyab \alpha n \rmyab
< \frac{m(\alpha-1)}{\alpha}+\frac{1}{\alpha}.
\]
Thus, by \eqref{ineq} again, we have
\[
f_\alpha (n)<1+\frac{m(\alpha-1)}{\alpha}+\frac{1}{\alpha}=m+1-\frac{m-1}{\alpha}<m+1
\]
and
\[
f_\alpha (n)>\frac{m(\alpha-1)}{\alpha}-\frac{1}{\alpha}> m-1,
\]
because $m<\alpha-1$. Therefore, $f_\alpha (n)=m$. \qedhere
\end{proof}
In particular, we get
\[
\Range(f_\varphi)=\{1\} \quad \text{and}\quad
\Range(f_\tau) = \bigl\llbracket 1 , 22 \bigr\rrbracket.
\]

A slight refinement of the first part of the proof of Proposition~\ref{golden} yields the following result.

\begin{proposition}
\label{p-golden}
Let $p\in\mathbb N$, $t\in\mathbb N$ and $\alpha(t)$ as in \eqref{alpha-special}. Then
\begin{equation}
\label{try-eq}
\Range\bigl(f_{p\mkern 1mu\alpha (t)}\bigr) \subseteq
\bigl\llbracket 1 , \lfloor p\mkern 1mu \alpha (t) \rfloor \bigr\rrbracket.
\end{equation}
\end{proposition}

\begin{proof}
Fix $p, t\in\mathbb N$ and simply write $\beta = p\mkern 1mu\alpha (t)$ for ease of notation. Let $n\in\mathbb N$. We begin by observing that $\beta^2-p\beta=p^2 t$, hence $\beta^2 n=p^2tn+p\beta n$. Therefore
\[
f_{\beta}(n)=p^2tn+\lfloor p\beta n\rfloor-\bigl\lfloor \beta  \lfloor \beta n \rfloor \bigr\rfloor.
\]
We will show the inclusion 
\(
\Range(f_{\beta}) \subseteq
\bigl\llbracket 1 , \lfloor \beta \rfloor \bigr\rrbracket.
\)
The main observation used in the estimates below is that
\begin{equation}
\label{obs}
p\lfloor x\rfloor\leq \lfloor px\rfloor\leq p\lfloor x\rfloor + (p-1)\quad \forall p\in\mathbb N,\quad \forall x\geq 0.
\end{equation}
Then, by \eqref{obs}, we have
\begin{align*}
f_\beta (n)
& < p^2tn+p\lfloor \beta n\rfloor + (p-1)-(\beta \lfloor\beta n\rfloor-1) \\
& =p^2tn+p-(\beta-p)\lfloor\beta n\rfloor\\
&=n(p^2t-\beta(\beta-p))+p+(\beta-p)\mkern 1mu \lmyab \beta n \rmyab \\
&=p+(\beta-p) \mkern 1mu \lmyab \beta n \rmyab \\
& < p+(\beta-p) \\
& = \beta.
\end{align*}
Similarly, using again \eqref{obs}, we can write
\begin{align*}
f_\beta (n)
& > p^2tn+p\lfloor\beta n\rfloor-\beta \lfloor\beta n\rfloor \\
& =p^2tn-(\beta-p)\lfloor\beta n\rfloor \\
&=p^2tn-(\beta-p)\bigl(\beta n- \lmyab \beta n \rmyab \bigr) \\
& =(\beta-p)\mkern 1mu \lmyab \beta n \rmyab \\ 
& > 0;
\end{align*}
the last strict inequality uses the fact that $\beta>p$ and $\beta\not\in\mathbb Q$. The required inclusion follows. 
\end{proof}

It is useful to note that in the preceding proof we obtained the stronger inequalities
\begin{equation} \label{beta-ineq}
(\beta-p)\mkern 1mu \lmyab \beta n \rmyab
< f_\beta (n)
< p+(\beta-p)\mkern 1mu \mkern 1mu \lmyab \beta n \rmyab, \quad \forall n\in \mathbb{N}.
\end{equation}

When we attempt to adapt the proof of Proposition~\ref{golden} to demonstrate the converse inclusion in \eqref{try-eq} for $p \geq 2$, the ``obvious'' choice of $\epsilon$ in Kronecker's approximation theorem, using now \eqref{beta-ineq} instead of \eqref{ineq}, significantly limits the range of $m$.

\subsection{Conjectures and partial results}

By Proposition~\ref{largest-range} we have
\begin{equation}  \label{pr-incl}
\Range(f_\alpha) \subseteq
\bigl\llbracket 0 , \lceil \alpha \rceil \bigr\rrbracket.
\end{equation}
Based on Wolfram Mathematica experimentation, see \cite{BC1}, we state the following

\begin{conjecture} \label{con}
For all positive irrational numbers $\alpha$, we have
\[
\bigl\llbracket 1 , \lfloor \alpha \rfloor \bigr\rrbracket \subseteq \Range(f_\alpha).
\]
\end{conjecture}

If Conjecture~\ref{con} is true, then we have the equality in Proposition~\ref{p-golden}.

Assuming Conjecture~\ref{con} holds, together with the inclusion in \eqref{pr-incl}, suggests the following quadruplicity of the range result.

\begin{proposition} \label{con2}
If Conjecture~{\rm\ref{con}} is true, then for all positive irrational parameters \(\alpha\), one of the following equalities holds:
\begin{enumerate}
\renewcommand*\theenumi{\Alph{enumi}}
\renewcommand*\labelenumi{\rm{(\theenumi)}}
  \item \label{cl-A}
\(\displaystyle
\Range(f_\alpha) = \bigl\llbracket 1 , \lfloor \alpha \rfloor \bigr\rrbracket\),
  \item \label{cl-B}
\(\displaystyle \Range(f_\alpha) =  \bigl\llbracket 0 , \lfloor \alpha \rfloor \bigr\rrbracket\),
  \item \label{cl-C}
\(\displaystyle
\Range(f_\alpha) =  \bigl\llbracket 1 , \lceil \alpha \rceil \bigr\rrbracket\),
  \item \label{cl-D}
\(\displaystyle
\Range(f_\alpha) = \bigl\llbracket 0 , \lceil \alpha \rceil \bigr\rrbracket\).
\end{enumerate}
\end{proposition}

\begin{proof}
The equality in (\ref{cl-A}) holds for all irrational parameters considered in Proposition~\ref{golden}. In particular, we have $\Range(f_{\varphi})=\{1\}$.

The equality in (\ref{cl-D}) holds for $\pi$ and $e$. We used Wolfram Mathematica  to calculate
\[
f_\pi(1) = 0, \ f_\pi(2) = f_\pi(3) = 1, \ f_\pi(4) = f_\pi(5) = 2, \ f_\pi(6) = 3, \ f_\pi(7) = 4,
\]
and
\begin{multline*}
f_e(7) = f_e(10) = 0, \ f_e(2) = f_e(3) = f_e(5) = f_e(6) = f_e(9) = 1, \\
f_e(1) = f_e(4) = f_e(8) = 2, \ f_e(11) = 3.
\end{multline*}

The equality in (\ref{cl-C}) holds for \(\alpha=\sqrt{2}\). Indeed,
\[
\Range\bigl(f_{\sqrt{2}}\bigr)=\{1,2\}.
\]
First, \(f_{\sqrt{2}}(k)=k\) for \(k\in\{1,2\}\), so \(\{1,2\}\subseteq\Range\bigl(f_{\sqrt{2}}\bigr)\).
Next, fix \(n\in\mathbb{N}\) and set \(m=\bigl\lfloor\sqrt{2}\,n\bigr\rfloor\).
By the definition of the floor, \(m\le \sqrt{2}\,n<m+1\). Since \(\sqrt{2}\,n\) is irrational, in fact \(m<\sqrt{2}\,n<m+1\).
Multiplying by \(\sqrt{2}\) gives
\[
2n-\sqrt{2}<\sqrt{2}\,m<2n,
\]
hence \(\bigl\lfloor\sqrt{2}\,m\bigr\rfloor\in\{2n-2,2n-1\}\). Therefore
\[
f_{\sqrt{2}}(n)
=2n-\bigl\lfloor \sqrt{2}\,\lfloor \sqrt{2}\,n\rfloor\bigr\rfloor
=2n-\bigl\lfloor \sqrt{2}\,m\bigr\rfloor
\in\{1,2\}.
\]
Thus \(\Range\bigl(f_{\sqrt{2}}\bigr)\subseteq\{1,2\}\), which, together with the converse inclusion, yields \(\Range\bigl(f_{\sqrt{2}}\bigr)=\{1,2\}\).

The equality in (\ref{cl-B}) holds for $\alpha = 2+\sqrt{2}$. That is,
\begin{equation*}  
\Range\bigl(f_{2+\sqrt{2}}\bigr) = \{0, 1, 2, 3\}.
\end{equation*}
The inclusion $\{0, 1, 2, 3\}\subseteq \Range\bigl(f_{2+\sqrt{2}}\bigr)$ follows from the equalities
\[
f_{2+\sqrt{2}}(1)=1, \ f_{2+\sqrt{2}}(2)=3, \ f_{2+\sqrt{2}}(3)=0, \ f_{2+\sqrt{2}}(4)=2.
\]

To prove the converse inclusion we use the following two properties of the floor function. For every positive real number \(x\) with
\(r = \bigl\lfloor 4 \bigl(x - \lfloor x \rfloor \bigr) \bigr\rfloor \in \{0,1,2,3\}\) we have
\begin{equation*}
  \lfloor 4 x \rfloor = 4 \lfloor x \rfloor + r.
\end{equation*}
For all positive integers \(n \in \mathbb{N}\) and all  positive real numbers \(x\) we have
\begin{equation*}
  n - \lfloor x \rfloor = \lceil n - x \rceil .
\end{equation*}

Let $n\in\mathbb N$ be arbitrary. A straightforward calculation shows that
\[
f_{2+\sqrt{2}}(n)
= 2n + \bigl\lfloor 4\sqrt{2}n \bigr\rfloor
- 2 \bigl\lfloor \sqrt{2}n \bigr\rfloor
- \bigl\lfloor 2\sqrt{2}n+\sqrt{2}\lfloor \sqrt{2}n\rfloor\bigr\rfloor.
\]
Set \(\lfloor \sqrt{2}n\rfloor = m\). Then there exists \(r \in \{0,1,2,3\}\) such that \(\lfloor 4\sqrt{2}n\rfloor = 4 m + r\). Hence,
\begin{align*}
f_{2+\sqrt{2}}(n)
& = 2 n + 4 m + r - 2 m - \bigl\lfloor 2\sqrt{2}n+\sqrt{2}m\bigr\rfloor \\
& = \Bigl\lceil 2 n + 2 m + r - 2\sqrt{2}n - \sqrt{2}m \Bigr\rceil \\
& = \lceil w \rceil.
\end{align*}

Next, we will obtain an estimate for the positive irrational number, we called it \(w\), whose ceiling is the value \(f_{2+\sqrt{2}}(n)\). First, it follows from \(\lfloor 4\sqrt{2}n\rfloor = 4 m + r\) that
\begin{align*}
m + \frac{r}{4} & < \sqrt{2} n < m + \frac{r+1}{4}, \\
\intertext{and, equivalently,}
\sqrt{2} n -  \frac{r+1}{4} & < m  < \sqrt{2} n - \frac{r}{4}.
\end{align*}
Transforming linearly the preceding two inequalities and adding them we obtain
\begin{equation*}
- \frac{r+1}{2} + \frac{r}{4}\sqrt{2}
< 2 n + 2m - 2\sqrt{2} n - \sqrt{2} m
< - \frac{r}{2} + \frac{r+1}{4}\sqrt{2}.
\end{equation*}
Adding \(r \in \{0,1,2,3\}\) to each term of the preceding inequalities and simplifying the left and righ expression, we get inequalities for \(w\)
\begin{equation*}
\frac{(2+\sqrt{2})r - 2}{4}
< w
< \frac{(2+\sqrt{2})r + \sqrt{2}}{4}.
\end{equation*}
Considering the extreme cases for \(r\in\{0,1,2,3\}\), that is \(r=0\) on the left-hand side and \(r=3\) on the right-hand side, we get the following open interval membership for \(w\):
\[
w \in \left(\mkern-3mu -\frac{1}{2}, \frac{3}{2} + \sqrt{2} \right).
\]
Consequently,
\[
f_{2+\sqrt{2}}(n) = \lceil w \rceil \in \{0,1,2,3\}.
\]
Since \(n\in\mathbb{N}\) was arbitrary, this proves $\Range\bigl(f_{2+\sqrt{2}}\bigr)\subseteq \{0, 1, 2, 3\}$, and completes the proof.
\end{proof}

The conjectured statements in Proposition~\ref{con2} account for the seemingly deep differences among positive  irrational numbers. It would be interesting to explore these differences further. We state two more conjectures.

\begin{conjecture} \label{con3}
For all irrational numbers $\alpha \in (0,1)$ we have
\[
\Range(f_\alpha) = \{0,1\}.
\]
That is, all irrational $\alpha \in (0,1)$ satisfy the equality in item {\rm (D)} in Proposition~{\rm\ref{con2}}.
\end{conjecture}

\begin{remark}
\label{rem23}
By Corollary~\ref{01range}, the preceding conjecture is equivalent to  $1\in\Range(f_\alpha)$, which in turn is equivalent to the following statement: For every \(x > 1\) there exists \(k \in \mathbb{N}\) such that
\[
\bigl\lfloor x \lceil x k \rceil \bigr \rfloor \geq \bigl\lceil x^2 k \bigr\rceil .
\]
\end{remark}

While Conjecture~\ref{con3} was formulated based on experiments in Wolfram Mathematica, we now aim to provide some partial rigorous evidence in support of its validity.

Our first observation is that Conjecture~\ref{con3} holds for all sub-unitary numbers $\alpha$ in some right-neighborhood of a square-root rational. We first set some notation. Let \(k,n \in \mathbb{N}\) be such that \(k < n\) and assume that \(\sqrt{k/n}\) is irrational. In this case, we call \(\sqrt{k/n}\) a \emph{square-root rational}, and denote it by $s_{k, n}$. With $k, n$ as above, we also write
\[
\epsilon_{k,n} = \frac{1}{n} \min\left\{
\bigl\lceil \sqrt{kn} \bigr\rceil - \sqrt{kn},
\sqrt{kn} - \bigl\lfloor \sqrt{kn} \bigr\rfloor  \right\} > 0.
\]

\begin{proposition}
\label{sqrtQ}
Let \(k,n \in \mathbb{N}\) be such that \(k < n\) and $s_{k, n}$ irrational. If $\alpha\in [s_{k, n}, s_{k, n}+\epsilon_{k, n})$, then $f_\alpha (n)=1$.
\end{proposition}

\begin{proof}
Let \(x, y \in \mathbb{R}\). The following three statements follow from the definitions of the floor and ceiling function:
\begin{equation*}
\lfloor x \rfloor < y \quad \Leftrightarrow  \quad
\lfloor x \rfloor < \lceil y \rceil \quad \Leftrightarrow \quad
 x  < \lceil y \rceil
\end{equation*}

Let \(k,n \in \mathbb{N}\) be such that \(k < n\) and assume that \(s_{k, n}=\sqrt{k/n}\) is irrational. Notice that \(s_{k, n}\) is irrational if and only if \(\sqrt{kn}\) is irrational. To evaluate \(f_{s_{k, n}}(n)\) we calculate
\[
\left\lfloor \left( \sqrt{\frac{k}{n}} \right)^2 n \right\rfloor = k,
\]
and estimate
\[
\Biggl\lfloor \sqrt{\frac{k}{n}} \biggl\lfloor \sqrt{\frac{k}{n}} \, n \biggr\rfloor \Biggr\rfloor
= \left\lfloor \sqrt{\frac{k}{n}} \bigl\lfloor \sqrt{k n} \bigr\rfloor \right\rfloor < \sqrt{\frac{k}{n}} \bigl\lfloor \sqrt{k n} \bigr\rfloor < \sqrt{\frac{k}{n}}   \sqrt{k n}  = k.
\]
Therefore \(f_{s_{k, n}}(n) = 1\).

Let \(x\) be a nonnegative real number such that $x<\epsilon_{k, n}$. In particular,
\[
x < \frac{1}{n} \Bigl(\bigl\lceil \sqrt{k n} \bigr\rceil - \sqrt{k n} \Bigr).
\]
Then
\[
\bigl\lfloor \sqrt{k n} \bigr\rfloor 
\leq \Biggl\lfloor\mkern -2mu \left(\sqrt{\frac{k}{n}} + x\right) n \Biggr\rfloor
\leq \sqrt{kn} + x n < \bigl\lceil \sqrt{k n} \bigr\rceil,
\]
and consequently,
\[
\bigl\lfloor s_{k, n}+x\bigr\rfloor 
= \Biggl\lfloor\mkern -2mu \left(\sqrt{\frac{k}{n}} + x\right) n \Biggr\rfloor =  \bigl\lfloor \sqrt{k n} \bigr\rfloor.
\]
Since we also have
\[
x< \frac{1}{n}\Bigl( \sqrt{k n} - \bigl\lfloor \sqrt{k n} \bigr\rfloor \Bigr)< \frac{1}{n} \, \frac{\sqrt{k n}}{\bigl\lfloor \sqrt{k n} \bigr\rfloor} \Bigl( \sqrt{k n} - \bigl\lfloor \sqrt{k n} \bigr\rfloor \Bigr),
\]
we can estimate
\begin{align*}
\Bigl\lfloor (s_{k, n}+x)\bigl\lfloor (s_{k, n}+x)n\bigr\rfloor\Bigr\rfloor 
& = \Biggl\lfloor\mkern -2mu \left(\sqrt{\frac{k}{n}} + x\right) \Biggl\lfloor \left(\sqrt{\frac{k}{n}} + x\right) n \Biggr\rfloor \Biggr\rfloor\\
& \leq  \left(\sqrt{\frac{k}{n}} + x\right) \bigl\lfloor \sqrt{k n} \bigr\rfloor   \\
& < \frac{1}{n} \left(
\sqrt{k n} + \frac{k n}{\bigl\lfloor \sqrt{k n} \bigr\rfloor} - \sqrt{k n}
\right) \bigl\lfloor \sqrt{k n} \bigr\rfloor  \\
& = k.
\end{align*}

Clearly,
\begin{equation} \label{eq-1st}
\bigl\lfloor (s_{k, n}+x)^2 n\bigr\rfloor=\Biggl\lfloor\mkern -2mu \left(\sqrt{\frac{k}{n}} + x\right)^2 n \Biggr\rfloor = k
\end{equation}
if and only if
\[
x < \frac{\sqrt{k+1} - \sqrt{k}}{\sqrt{n}}
= \frac{1}{\sqrt{n}\bigl(\sqrt{k+1} + \sqrt{k}\bigr)}.
\]
Since
\[
x <\epsilon_{k, n} < \frac{1}{2n}< \frac{1}{\sqrt{n}\bigl(\sqrt{k+1} + \sqrt{k}\bigr)},
\]
we deduce that \eqref{eq-1st} holds with so chosen \(x\). Since under this estimate we also have
\[
\Bigl\lfloor (s_{k, n}+x)\bigl\lfloor (s_{k, n}+x)n\bigr\rfloor\Bigr\rfloor<k,
\]
we have proved that for
\(
\alpha = s_{k, n} + x
\)
we have
\(
f_\alpha(n) = 1.
\)
\end{proof}

Our second observation in support of the validity of Conjecture~\ref{con3} is the calculation of the sets
\[
A_m = \bigl\{ \alpha \in (0,1] :  f_\alpha(m) = 1 \bigr\},
\]
where \(m \in \mathbb{N}\setminus\{1\}\). With this notation, and using Proposition~\ref{subunitary}, Conjecture~\ref{con3}  can be restated as
\[
\bigcup \bigl\{A_m  : m\in\mathbb{N}\setminus\{1\} \bigr\}= (0,1) \setminus \Bigl\{ \frac{1}{b} : b \in \mathbb{N}  \Bigr\}.
\]

We wrote Wolfram Mathematica code that calculates the sets \(A_m\) and their finite unions. For example,
\begin{align*}
\bigcup_{m=2}^{50} A_m & =
\mkern 160mu \left[ \frac{1}{5 \sqrt{2}} , \frac{1}{7} \right) \\
&\mkern 50mu \medcupd
\left[ \frac{1}{4 \sqrt{3}} , \frac{7}{48} \right) \medcupd
\left[ \frac{1}{\sqrt{47}} , \frac{1}{6} \right) \\
&\mkern 50mu \medcupd  \left[ \frac{1}{\sqrt{35}} , \frac{6}{35} \right) \medcupd
\left[ \frac{1}{\sqrt{34}} , \frac{1}{5} \right) \\
&\mkern 50mu \medcupd
\left[ \frac{\sqrt{2}}{7} , \frac{10}{49}\right)
\medcupd \left[ \frac{1}{2 \sqrt{6}} , \frac{1}{4} \right) \\
&\mkern 50mu \medcupd
\left[ \sqrt{\frac{3}{47}} , \frac{1}{3}\right) \medcupd
\left[  \sqrt{\frac{5}{44}}  , \frac{1}{2} \right) \medcupd
\left[ 2 \sqrt{\frac{3}{47}} , 1 \right) \\
& \approx \mkern 208mu \bigl[0.14142, 0.14286\bigr) \\
&\mkern 50mu \cup \bigl[0.14434, 0.14583\bigr) \cup
\bigl[0.14586, 0.16667\bigr) \\
&\mkern 50mu \cup \bigl[0.16903, 0.17143 \bigr) \cup \bigl[0.1715, 0.2\bigr) \\
&\mkern 50mu \cup \bigl[0.20203, 0.20408 \bigr) \cup \bigl[0.20412, 0.25 \bigr) \\
&\mkern 50mu \cup \bigl[0.25265, 0.33333 \bigr)
\cup \bigl[0.3371, 0.5 \bigr) \cup \bigl[0.50529,1\bigr)
\end{align*}

\begin{figure}[htbp]
  \centering
  \includegraphics[width=\textwidth]{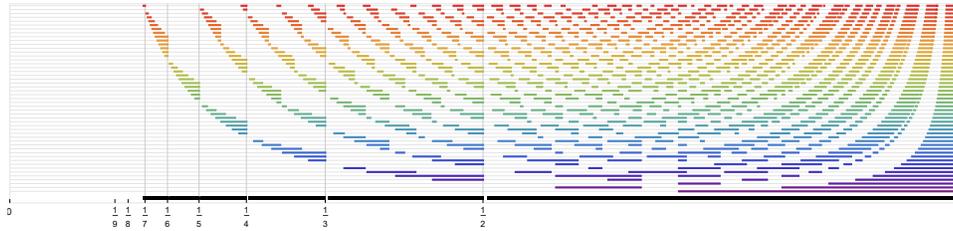}
  \caption{The sets \(A_2,\ldots,A_{50}\) in rainbow colors and their union in black}
  \label{fig:union50}
\end{figure}

Figure~\ref{fig:union50} displays the individual sets \(A_m\) for \(m \in \{2,\ldots,50\}\), rendered in rainbow colors and arranged vertically, with \(A_2\) at the bottom and \(A_{50}\) at the top. The union of these sets is shown in black at the very bottom of the image. The union consists of \(10\) intervals. Among the \(9\) gaps between successive intervals, three are particularly narrow. These occur near \(0.1458\), \(0.1715\), and \(0.2041\), with respective widths of approximately \(0.00003\), \(0.00007\), and \(0.00004\). Even at a magnification of 5000\%, only the second gap is barely discernible in the PDF image.

The only positive irrational number \(\alpha\) for which the function \(f_\alpha\) is constant that we encountered in our explorations is \(\alpha = \varphi\), the golden ratio. Therefore we conjecture:

\begin{conjecture} \label{con4}
For all positive irrational $\alpha$ we have
\[
\Range(f_\alpha) = \{1\} \qquad \Leftrightarrow \qquad \alpha = \varphi.
\]
\end{conjecture}

\begin{remark}
We disprove the following invariance conjecture: For all irrational \(\alpha > 1\) and all \(p \in \mathbb{N}\), if \(\Range(f_{\alpha})\) satisfies the equality in item {\rm(X)}, where {\rm X} $\in$ $\{${\rm A, B, C, D}$\}$, then \(\Range(f_{p \alpha})\) would also satisfy the same equality.

To disprove this conjecture, consider $\alpha = 1+1/\sqrt{2}$ and calculate
\[
f_\alpha(1) = 1, \quad f_\alpha(2) = 0, \quad f_\alpha(7) = 2.
\]
Since \(\lceil\alpha\rceil = 2\), by Proposition~\ref{largest-range} we have \(\Range(f_\alpha) = \{0,1,2\}\). Hence, $\alpha = 1+1/\sqrt{2}$ satisfies the equality in item (\ref{cl-D}) in Proposition~{\rm\ref{con2}}. Notice that \(2 \alpha = 2+\sqrt{2}\). We proved in the proof of Proposition~\ref{con2} that \(2\alpha = 2+\sqrt{2}\) satisfies the equality in item (\ref{cl-B}) in Proposition~{\rm\ref{con2}}.
\end{remark}

We conclude with a couple of remarks regarding variations of the function $f_\alpha$ discussed here.

\begin{remark}
Throughout this article, we worked with the domain of the function $f_\alpha$ being the set of positive integers. Naturally, one may ask what happens if we consider the domain as being the set of negative integers. For this, we note that, since $\lceil x\rceil=-\lfloor -x\rfloor$, given $n\in\mathbb N$ we have
\[
f_\alpha(-n)=-g_\alpha(n),
\]
where now $g_\alpha: \mathbb N\to\mathbb Z,\, g_\alpha (n)=\bigl\lceil \alpha^2 \, n \bigr\rceil - \bigl\lceil \alpha  \lceil \alpha n \rceil \bigr\rceil.$ In other words, $f_\alpha(-\mathbb N)=-g_\alpha(\mathbb N)$; the analysis of the $\Range(g_\alpha)$ will most likely be of a similar flavor to one we have seen for $\Range(f_\alpha)$.
\end{remark}

\begin{remark}
Finally, one may ask about negative parameters $\alpha$. Note now that for $\alpha<0$, we can write \(f_{\alpha}=h_{-\alpha}\) where, for $n\in\mathbb N$ and $\beta>0,$ we have denoted
\[
h_{\beta}(n)=\bigl\lfloor \beta^2 \, n \bigr\rfloor - \bigl\lfloor \beta  \lceil \beta n \rceil \bigr\rfloor.
\]
This is yet another variation on the function $f_\alpha$. There are a few other such variants, obtained by appropriately interchanging the floor and ceiling operations in the definition of $f_\alpha$. The analysis of their ranges is left to the interested reader.
\end{remark}

\section*{Acknowledgement}
\'A.B.  acknowledges the support from an  AMS-Simons Research Enhancement Grant for PUI Faculty.

\end{document}